# UNIQUENESS THRESHOLDS ON TREES VERSUS GRAPHS

By Allan Sly

*University of California, Berkeley*

Counter to the general notion that the regular tree is the worst case for decay of correlation between sets and nodes, we produce an example of a multi-spin interacting system which has uniqueness on the $d$-regular tree but does not have uniqueness on some infinite $d$-regular graphs.

## 1. Definitions.

DEFINITION 1. On a graph $G = (V, E)$ a *Gibbs measure* (also *Markov random field* or *graphical model*) is a distribution $\sigma$ taking values in $\mathcal{C}^V$, for some finite set $\mathcal{C}$, which satisfies the Markov property

$$P(\sigma_A, \sigma_B | \sigma_S) = P(\sigma_A | \sigma_S) P(\sigma_B | \sigma_S),$$

when $A, B$ and $S$ are disjoint subsets of $V$ such that every path in $G$ from $A$ to $B$ passes through $S$ and where $\sigma_U$ denotes the natural projection of $\sigma$ from $\mathcal{C}^V$ to $\mathcal{C}^U$ for $U \subset V$.

Throughout this paper we restrict our attention to Gibbs measure with fixed activities and interactions on the vertices and edges so that

$$P(\sigma) = \frac{1}{Z} \exp\left[\sum_{v \in V} g(\sigma_v) + \sum_{(u,v) \in E} h(\sigma_u, \sigma_v)\right]$$

where $Z$ is a normalizing constant and $g : \mathcal{C} \to \mathbb{R}$ and $h : \mathcal{C}^2 \to \mathbb{R} \cup \{-\infty\}$ are functions.

DEFINITION 2. On an infinite graph $G = (V, E)$ with finite degrees we say that a Gibbs measure has *uniqueness* if for any finite set $A \subset V$,

$$\lim_n \sup_{\sigma^1, \sigma^2} d_{TV}(P(\sigma_A = \cdot | \sigma_{S(A,n)} = \sigma^1), P(\sigma_A = \cdot | \sigma_{S(A,n)} = \sigma^2)) = 0,$$

where $d_{TV}$ is the total variation distance, $S(A, n)$ denotes the set $\{u \in T : d(u, A) = n\}$ and the supremum is taken over all boundary conditions $\sigma^1, \sigma^2 \in \mathcal{C}^{S(v,A)}$ on $S(v, A)$.

**2. Background.** Sokal [9] conjectured that uniqueness in the hard-core model on the $d$-regular tree implies uniqueness on any $d$-regular graph. He also speculated that this might also be true for random colorings. Mossel [7] suggested that this may in fact hold for every spin system.

Determining the regimes for uniqueness and nonuniqueness on regular trees can often be done through recursions and so can be easier than on general graphs. For many systems there exists a threshold in the parametrization between the uniqueness and nonuniqueness regimes. The correctness of the conjectures would then allow us to deduce uniqueness of Gibbs measures for regular graphs in the regime for which uniqueness holds on the regular tree. In general, determining the thresholds for uniqueness in regular graphs is hard.

The intuition behind such conjectures is that the regular tree has the most vertices at distance $n$ from the root and so this boundary has the greatest influence on the root. In this sense loops constitute wasted influence. However, loops in the graph create extra dependence between the states of the neighbors. This is crucial in the construction of our counterexample.

Weitz in [10] showed that marginals of the hard-core model on a $d$-regular graph could be exactly evaluated by calculating the marginals on a tree of self-avoiding random walks. This approach establishes efficient deterministic polynomial time algorithms for approximately counting independent sets on $d$-regular graphs. In [10] the generalization to general 2-spin systems is implicitly given, though [5] gives an explicit description. Tree-based constructions for spin systems have also been used in [1, 4] and [8].

As an immediate consequence of this construction [10] shows that any 2-spin system which has strong spatial mixing on the $d$-regular tree also has strong spatial mixing on all graphs of maximum degree $d$. That is, the worst case for strong spatial mixing is the $d$-regular tree. In the hard-core and antiferromagnetic Ising models [10] also showed that uniqueness on the $d$-regular tree in fact implies strong spatial mixing on the $d$-regular tree and so implies uniqueness on all graphs of maximum degree $d$, proving the conjecture of Sokal. This approach does not apply more generally to all 2-spin systems; in particular in the ferromagnetic Ising model uniqueness on the $d$-regular tree does not in general imply strong spatial mixing (see, e.g., [10]). It is unknown whether in 2-spin systems uniqueness on the $d$-regular tree implies uniqueness on all $d$-regular graphs.

Our counterexample leaves open two important related conjectures:

- It remains open whether strong spatial mixing on the $d$-regular tree implies strong spatial mixing on all graphs of maximum degree $d$.

- In the specific example of random colorings it remains an open problem to show that $q \geq d+1$ colors implies uniqueness on all graphs with maximum degree $d$.

Later we will discuss why our example does not provide a counterexample to the first conjecture. It remains of interest to show when uniqueness on the $d$-regular tree implies uniqueness on all graphs of maximum degree $d$, for instance whether this holds on monotone spin systems.

Determining regimes for uniqueness and strong spatial mixing plays a key role in analyzing the performance of algorithms for sampling from Markov chains and approximately counting distributions. On lattices it is known that strong spatial mixing implies rapid mixing of the Glauber dynamics [3, 6]. As determining the thresholds for strong spatial mixing on trees is in many instances simpler than on the lattice, such a result would be a powerful tool for determining on which graphs mixing is rapid. For instance, in the 2-spin setting, Weitz [10] was able to improve the best known bound for rapid mixing in the hard-core model on $\mathbb{Z}^2$ from $\lambda < 1.508$ to $\lambda < 1.6875$.

We should note that other seemingly counterintuitive uniqueness results have been found in other settings. For example, it was shown in [2] that for some graphs uniqueness in the hard-core model need not be monotone in the activity parameter. In fact it was shown that nonuniqueness on a subgraph need not imply non-uniqueness on the whole graph, even on trees, and so uniqueness is not monotone in the degrees of the graph.

**3. Our construction.** The introduction of multiple spins adds greater complexity to the question of spatial mixing of a graph. In a 2-spin system on a tree the marginal at the root is maximized by maximizing the marginal of some state in each of its neighbors. When multiple spins are involved the whole collection of spins of the neighbors determines the marginal at the root. For instance, for random colorings the marginal at the root is determined by which colors do not appear amongst its neighbors. In a tree the colors of the neighbors are conditionally independent given the color at the root. By contrast, in a graph the colors at the neighbors can be conditionally dependent increasing the total number of distinct colors. This observation is the source of motivation for our model. In the graph the dependence between the neighbors of the root more than makes up for the smaller number of vertices on the boundary.

Our model begins with an antiferromagnetic Potts model taking states $\{1, \ldots, q\}$,

$$P(\sigma) = \frac{1}{Z} \exp\left[-\beta \sum_{(u,v) \in E} 1_{\{\sigma_u = \sigma_v\}}\right]$$

on the $q$-ary tree $T$ with root $v$. A simple modification of Lemma 3 below shows that for large $q$ and any $0 < \beta < \infty$,

$$\limsup_n \sup_{\sigma^*} P(\sigma_v = 1 | \sigma_{S(v,n)} = \sigma^*) \leq \tfrac{1}{4} \tag{1}$$

where $S(v,n)$ denotes $\{u \in T : d(u,v) = n\}$ and $\sigma^*$ is any boundary condition on $S(v,n)$. In fact we can replace $1/4$ with $C/q$ for some constant $C$. By contrast one can show that on the $q+1$ regular graph $G$ described in Section 4

$$\sup_{0 < \beta < \infty} \inf_n \sup_{\sigma^*} P(\sigma_v = 1 | \sigma_{S(v,n)} = \sigma^*) = 1 \tag{2}$$

since it is constructed so that the neighbors of $v$ coordinate their spins to take different values.

This antiferromagnetic Potts model does not have uniqueness on the $(q+1)$-regular tree for $\beta > \ln(q+1)$. However, we exploit the difference in (1) and (2) by adding an extra state with different interactions and intensity so that the modified system has uniqueness on the $q+1$-regular tree but not on the graph $G$. This allows us to prove the main result.

**4. Main result.** Let $G = (V, E)$ be a finite graph. Denote the possible states as elements of $[q+1] = \{1, 2, \ldots, q+1\}$ and let

$$P(\sigma) = \frac{1}{Z} \exp\left[ \lambda \sum_{u \in V} 1_{\{\sigma_u = q+1\}} - \beta \sum_{(u,v) \in E} \eta(\sigma_u, \sigma_v) \right],$$

where $\lambda$ and $\beta$ are positive and

$$\eta(i,j) = \begin{cases} 1, & \text{if } 1 \leq i = j \leq q, \\ 1, & \text{if } i = q+1, j \in \{1,2\}, \\ 1, & \text{if } j = q+1, i \in \{1,2\}, \\ 0, & \text{otherwise.} \end{cases}$$

The first $q$ states are just a standard antiferromagnetic Potts model.

THEOREM 1. *For $q \geq 90$ and $e^\lambda > q^3 4^{q+1}$ and for any $0 < \beta < \infty$, the Gibbs measure is unique on the $q$-ary tree.*

Now consider the same model on an infinite graph $G$ defined recursively as follows:

- Start with one root vertex $v$ in row 1.
- For each vertex $u$ in row $i$, when $i$ is odd add a $q-1$-clique of vertices $u_1, \ldots, u_{q-1}$ in row $i+1$ and connect each of them to $u$ to form a $q$-clique.
- For each $q-1$-clique $u_1, \ldots, u_{q-1}$ in row $i$, when $i$ is even add $q-1$ vertices $w_1, \ldots, w_{q-1}$ in row $i+1$ and connect $u_l$ to $w_l$ and $w_{l+1}$ for $1 \leq l \leq q-2$ and connect $u_{q-1}$ to $w_1$ and $w_{q-1}$.

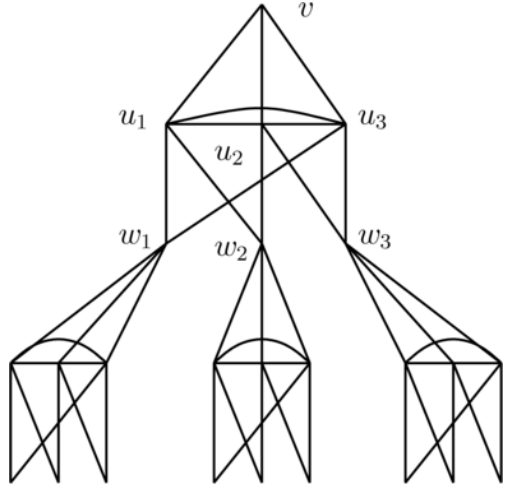

Fig. 1. *The first five rows of $G$ with $q = 4$.*

The graph $G$ is $q+1$-regular except at $v$. See Figure 1.

THEOREM 2. *For sufficiently large finite $\beta$, uniqueness does not hold on the graph $G$.*

These two results show that by taking $\beta$ large enough we have an example with uniqueness on the $d$-regular tree but not on the $d$-regular graph $G$.

**5. Proofs.** We prove four lemmas about the Gibbs measure on the $q$-ary tree and conclude with the proof of Theorem 1. We will denote the $q$-ary tree as $T$ with root $v$, which has children $u_1, \ldots, u_q$.

LEMMA 1. *For $1 \leq l \leq q$ let $T_l$ denote the subgraph generated by $u_l$ and all its descendants and let*

$$P^l(i) = P^{T_l}(\sigma_{u_l} = i | \sigma_{S(v,n) \cap T_l} = \sigma^*)$$

*where $P^{T_l}$ is the restriction to the subgraph $T_l$ disconnected from the rest of the tree and $\sigma^*$ is any boundary condition. Then for any $0 < \beta < \infty$ and $1 \leq k \leq q$,*

(3)
$$P(\sigma_v = k | \sigma_{S(v,n)} = \sigma^*)$$
$$\leq \frac{1}{1 + \sum_{j \in \{3,\ldots,q\}, j \neq k} \prod_{l=1}^{q}(1 - P^l(j))}.$$

Proof. In our model we have

$$P(\sigma_v = k | \sigma_{S(v,n)} = \sigma^*) = \frac{\sum_{j' \in [q+1]^q} \psi(k, j') \prod_{l=1}^q P^l(j'_l)}{\sum_{j \in [q+1]} \sum_{j' \in [q+1]^q} \psi(j, j') \prod_{l=1}^q P^l(j'_l)}, \quad (4)$$

where $j' = (j'_1, \ldots, j'_q) \in [q+1]^q$ represents the states of $u_1, \ldots, u_q$ and where $\psi(j, j')$ is given by

$$\exp\left[\lambda 1_{\{j=q+1\}} - \beta \sum_{l=1}^q \eta(j, j'_l)\right]. \quad (5)$$

Substituting this expression for $\psi$, we get

$$\sum_{j' \in [q+1]^q} \psi(j, j') \prod_{l=1}^q P^l(j'_l)$$

$$= \exp(\lambda 1_{\{j=q+1\}}) \sum_{j' \in [q+1]^q} \prod_{l=1}^q \exp(-\beta \eta(j, j_l)) P^l(j'_l)$$

$$= \exp(\lambda 1_{\{j=q+1\}}) \prod_{l=1}^q \sum_{j_l \in [q+1]} \exp(-\beta \eta(j, j_l)) P^l(j'_l)$$

$$(6) \quad = \begin{cases} \prod_{l=1}^q (1 - (1 - e^{-\beta})(P^l(j) + P^l(q+1))), & \text{if } j \in \{1, 2\}, \\ \prod_{l=1}^q (1 - (1 - e^{-\beta}) P^l(j)), & \text{if } 3 \leq j \leq q, \\ e^\lambda \prod_{l=1}^q (1 - (1 - e^{-\beta})(P^l(1) + P^l(2))), & \text{if } j = q+1. \end{cases}$$

It follows that for $3 \leq j \leq q$,

$$\sum_{j' \in [q+1]^q} \psi(j, j') \prod_{l=1}^q P^l(j'_l) \geq \prod_{l=1}^q (1 - P^l(j)) \quad (7)$$

and for $1 \leq j \leq q$,

$$\sum_{j' \in [q+1]^q} \psi(j, j') \prod_{l=1}^q P^l(j'_l) \leq \sum_{j' \in [q+1]^q} \prod_{l=1}^q P^l(j'_l) = 1. \quad (8)$$

Applying (7) and (8) to (4) we get the inequality

$$P(\sigma_v = k | \sigma_{S(v,n)} = \sigma^*)$$

$$\leq \frac{1}{1 + \sum_{j \in \{3, \ldots, q\}, j \neq k} \prod_{l=1}^q (1 - P^l(j))}. \quad \square$$

LEMMA 2. *Assume we have that for all $1 \leq j \leq q$ and $1 \leq l \leq q$, $P^l(j) \leq p$ for some $\frac{1}{2} \leq p < 1$. Then if $0 < \beta < \infty$, for any $1 \leq k \leq q$,*

$$\sum_{j \in \{3,\ldots,q\}, j \neq k} \prod_{l=1}^{q} (1 - P^l(j)) \geq (q/6 - 3)p(1-p). \tag{9}$$

PROOF. Suppose that we fix $P^l$ for $1 \leq l \leq q-1$ and suppose we want to minimize

$$\sum_{j \in \{3,\ldots,q\}, j \neq k} \prod_{l=1}^{q} (1 - P^l(j)) = - \sum_{j \in \{3,\ldots,q\}, j \neq k} P^q(j) \prod_{l=1}^{q-1} (1 - P^l(j)) \tag{10}$$

$$+ \sum_{j \in \{3,\ldots,q\}, j \neq k} \prod_{l=1}^{q-1} (1 - P^l(j)) \tag{11}$$

in $P^q$. This is a linear equation in the $P^q(j)$. For $1 \leq j \leq q$ the coefficient of $P^q(j)$ is $-\prod_{l=1}^{q-1}(1 - P^l(j)) \leq 0$ while the coefficient of $P^q(q+1)$ is 0. So the minimum, subject to the constraints in the hypothesis, must have $P^q(q+1) = 0$ and for some states $1 \leq j_1, j_2 \leq q$, $P^q(j_1) = 1 - P^q(j_2) = p$. This can be applied to $P^l$ for any $1 \leq l \leq q$ so the left-hand side of (9) is minimized by taking

$$P^l(j_{l1}) = 1 - P^l(j_{l2}) = p \tag{12}$$

for some choice of states $1 \leq j_{11}, \ldots, j_{q1}, j_{12}, \ldots, j_{q2} \leq q$.

So assume the $P^l$ are of the form given in (12). Let $A$ be the set of states in $\{j \in \{3, 4, \ldots, q\}, j \neq k\}$ that appear at most once in the list $j_{11}, \ldots, j_{q1}$ and at most twice in the combined list $j_{11}, \ldots, j_{q1}, j_{12}, \ldots, j_{q2}$. If $j \in A$, then $P^l(j) = 0$ for all but at most two values of $l$ so the product $\prod_{l=1}^{q}(1 - P^l(j))$ has at most two terms not equal to 1. The nonone terms in the product are either equal to $p$ or $1-p$ and at most one of them is $1-p$. Then since $1 > p \geq 1 - p$ we have

$$\prod_{l=1}^{q} (1 - P^l(j)) \geq p(1-p)$$

for any $j \in A$. The proof will be completed by showing that $|A| \geq q/6 - 3$ so that

$$\sum_{j \in \{3,\ldots,q\}, j \neq k} \prod_{l=1}^{q} (1 - P^l(j)) \geq |A|p(1-p) \geq (q/6-3)p(1-p).$$

Let $a_1$ (resp. $a_2$) be the number of states that appear exactly once (resp. at least twice) in $j_{11}, \ldots, j_{q1}$. Let $b_2$ (resp. $b_3$) be the number of states that appear exactly twice (resp. at least three times) in $j_{12}, \ldots, j_{q2}$. By counting states according to how many times they appear in $j_{11}, \ldots, j_{q1}$ and $j_{12}, \ldots, j_{q2}$ we have

(13)  $$a_1 + 2a_2 \leq q, \qquad 2b_2 + 3b_3 \leq q.$$

Any state not in $|A|$ must appear either at least twice in $j_{11}, \ldots, j_{q1}$ or at least three times in $j_{12}, \ldots, j_{q2}$ or once in $j_{11}, \ldots, j_{q1}$ and twice in $j_{12}, \ldots, j_{q2}$ and so

$$\begin{aligned}
|A| &\geq q - 3 - [a_2 + b_3 + \min\{a_1, b_2\}] \\
&\geq q - 3 - \left[\frac{q - a_1}{2} + \frac{q - 2b_2}{3} + \min\{a_1, b_2\}\right] \\
&= q/6 - 3 - [\min\{a_1, b_2\} - a_1/2 - 2b_2/3] \\
&\geq q/6 - 3.
\end{aligned}$$

The second inequality follows from (13) while the final inequality follows from the fact that $\min\{a_1, b_2\} \leq (a_1 + b_2)/2 \leq a_1/2 + 2b_2/3$. □

LEMMA 3. *If $q \geq 90$ and $e^\lambda > q2^q$, then for any $0 < \beta < \infty$ there exists an $N = N(q, \lambda, \beta)$ such that for all $n > N$,*

$$\inf_{\sigma^*} P(\sigma_v = q + 1 | \sigma_{S(v,n)} = \sigma^*) \geq \tfrac{1}{2},$$

*where the infimum is over all boundary conditions on $S(v,n)$.*

PROOF. Define

$$p_n = \sup_{\sigma^*} \sup_{1 \leq i \leq q} P(\sigma_v = i | \sigma_{S(v,n)} = \sigma^*),$$

where the supremum is taken over all boundary conditions $\sigma^*$ on $S(v, n)$. Then by Lemmas 1 and 2 for $1 \leq k \leq q$ we have $P(\sigma_v = k | \sigma_{S(v,n)} = \sigma^*) \leq f(p_{n-1})$ where $f(p) = 1/(1 + (q/6 - 3)p(1 - p))$ and so $p_n \leq f(p_{n-1})$. Now $p = f(p)$ has three solutions, $1, 1/\sqrt{q/6 - 3}, -1/\sqrt{q/6 - 3}$. Since $\beta < \infty$ the constraints are soft and $0 < p_1 < 1$. Since

$$\left.\frac{d}{dp} f(p)\right|_{p=1} = q/6 - 3 \geq 12,$$

then $p_n$ must decrease toward $1/\sqrt{q/6 - 3}$, which is less than $\tfrac{1}{2}$, until for some $n$, $p_n < \tfrac{1}{2}$. Now when $p_{n-1} < \tfrac{1}{2}$ we have that $p_n \leq f(\tfrac{1}{2})$ so for large enough $n$ we have that

$$p_n \leq \frac{1}{1 + (q/6 - 3)/4} \leq \frac{1}{4}.$$

Now given that $P^l(1), P^l(2) \leq \frac{1}{4}$ by applying (6) and (8) to (4) we get

$$P(\sigma_v = q+1 | \sigma_{S(v,n)} = \sigma^*) \geq \frac{e^\lambda \prod_{l=1}^q (1 - P^l(1) - P^l(2))}{q + e^\lambda \prod_{l=1}^q (1 - P^l(1) - P^l(2))}$$

$$\geq \frac{2^{-q} e^\lambda}{q + 2^{-q} e^\lambda}$$

$$\geq \frac{1}{2}$$

as required. □

Lemma 3 established that there are eventually a large proportion of vertices in state $q + 1$ when sufficiently far from the boundary. This enables us to establish a contraction mapping. Now suppose we consider

$$P^l(j) = P^{T_l}(\sigma_{u_l} = j | \sigma_{S(v,n)} = \sigma^*)$$

as an element of $\mathbb{R}^{q+1}$ and let $P = P(j) \in \mathbb{R}^{q+1}$ denote $P(\sigma_v = j | \sigma_{S(v,n)} = \sigma^*)$, the marginal at $v$. Then $P$ can be derived from $P^1, \ldots, P^q$ by $P = g(P^1, \ldots, P^q)$ where

$$g(P^1, \ldots, P^q)(j) = \frac{\sum_{j' \in [q+1]^q} \psi(j, j') \prod_{l=1}^q P^l(j'_l)}{\sum_{j \in [q+1]} \sum_{j' \in [q+1]^q} \psi(j, j') \prod_{l=1}^q P^l(j'_l)}.$$

LEMMA 4. *Let $q \geq 90$, $e^\lambda > q^3 4^{q+1}$ and $0 < \beta < \infty$. Suppose we have distributions $P^1, \ldots, P^q$ and $Q^1, \ldots, Q^q$ all satisfying $P^l(q+1), Q^l(q+1) \geq \frac{1}{2}$. If $P = g(P^1, \ldots, P^q)$ and $Q = g(Q^1, \ldots, Q^q)$, then*

$$\|P - Q\|_1 \leq C \max_{1 \leq l \leq q} \|P^l - Q^l\|_1,$$

*where $0 < C < 1$ and $\|\cdot\|_1$ is the usual $L^1$ norm on $\mathbb{R}^{q+1}$ (or equivalently in the total variation distance between the distributions).*

PROOF. Denote $K_P$ as

$$K_P = \sum_{j \in [q+1]} \sum_{j' \in [q+1]^q} \psi(j, j') \prod_{l=1}^q P^l(j'_l)$$

and denote $K_Q$ similarly. Now observe the simple inequality that if $0 \leq x_1, \ldots, x_q \leq 1$ and $0 \leq y_1, \ldots, y_q \leq 1$, then

(14)
$$\left| \prod_{l=1}^q x_l - \prod_{l=1}^q y_l \right| = \left| \sum_{j=1}^q (x_j - y_j) \prod_{l=1}^{j-1} x_l \prod_{l=j+1}^q y_l \right|$$

$$\leq \sum_{j=1}^q |x_j - y_j|.$$

Applying (14) to (6), it follows that for $1 \leq j \leq q$,

$$\left| \sum_{j' \in [q+1]^q} \psi(j,j') \prod_{l=1}^{q} P^l(j'_l) - \sum_{j' \in [q+1]^q} \psi(j,j') \prod_{l=1}^{q} Q^l(j'_l) \right| \leq \sum_{l=1}^{q} \|P^l - Q^l\|_1$$

and that

$$\left| \sum_{j' \in [q+1]^q} \psi(q+1,j') \prod_{l=1}^{q} P^l(j'_l) - \sum_{j' \in [q+1]^q} \psi(q+1,j') \prod_{l=1}^{q} Q^l(j'_l) \right|$$

$$\leq e^\lambda \sum_{l=1}^{q} \|P^l - Q^l\|_1$$

and so

$$|K_P - K_Q| \leq (q + e^\lambda) \sum_{l=1}^{q} \|P^l - Q^l\|_1.$$

Also note that

$$K_P, K_Q \geq e^\lambda \prod_{l=1}^{q} P^l(q+1) \geq e^\lambda 2^{-q}.$$

Then for $1 \leq j \leq q$ using these estimates,

$$|P(j) - Q(j)|$$

$$= \left| \frac{\sum_{j' \in [q+1]^q} \psi(j,j') \prod_{l=1}^{q} P^l(j'_l)}{K_P} - \frac{\sum_{j' \in [q+1]^q} \psi(j,j') \prod_{l=1}^{q} Q^l(j'_l)}{K_Q} \right|$$

$$= \left| \frac{\sum_{j' \in [q+1]^q} \psi(j,j') \prod_{l=1}^{q} P^l(j'_l) - \sum_{j' \in [q+1]^q} \psi(j,j') \prod_{l=1}^{q} Q^l(j'_l)}{K_P} \right.$$

$$\left. + \frac{(K_Q - K_P) \sum_{j' \in [q+1]^q} \psi(j,j') \prod_{l=1}^{q} Q^l(j'_l)}{K_P K_Q} \right|$$

$$\leq \frac{\sum_{l=1}^{q} \|P^l - Q^l\|_1}{e^\lambda 2^{-q}} + \frac{(q + e^\lambda) \sum_{l=1}^{q} \|P^l - Q^l\|_1}{(e^\lambda 2^{-q})^2}$$

$$\leq e^{-\lambda}(2^q + 4^q(1 + qe^{-\lambda})) \sum_{l=1}^{q} \|P^l - Q^l\|_1.$$

Since $P(q+1) = 1 - \sum_{j=1}^{q} P(j)$ we have

$$\|P - Q\|_1 \leq 2 \sum_{j=1}^{q} |P(j) - Q(j)|$$

$$\leq 2qe^{-\lambda}(2^q + 4^q(1+qe^{-\lambda}))\sum_{l=1}^{q}\|P^l - Q^l\|_1$$

$$\leq 2q^2 e^{-\lambda}(2^q + 4^q(1+qe^{-\lambda}))\max_{1\leq l\leq q}\|P^l - Q^l\|_1,$$

which establishes the result since $2q^2 e^{-\lambda}(2^q + 4^q(1+qe^{-\lambda})) < e^{-\lambda}q^3 4^{q+1} < 1$.
□

PROOF OF THEOREM 1. Combining Lemmas 3 and 4 shows that for any $1 \leq i \leq q+1$ and $\varepsilon > 0$, then for large enough $n$,

$$\sup_{\sigma^{*1}, \sigma^{*2}} |P(\sigma_v = i|\sigma_{S(v,n)} = \sigma^{*1}) - P(\sigma_v = i|\sigma_{S(v,n)} = \sigma^{*2})| < \varepsilon$$

which is sufficient to establish uniqueness and prove Theorem 1. □

**6. Proof of Theorem 2.** We will show that

(15) $$\sup_{0<\beta<\infty} \inf_n P(\sigma_v = 1|\sigma_{S(v,2n)} \equiv 1) = 1.$$

Fix $0 < \beta < \infty$ and let

$$p_n = P(\sigma_v = 1|\sigma_{S(v,2n)} \equiv 1).$$

Let $u_1, \ldots, u_{q-1}$ denote the vertices in the second row and let $w_1, \ldots, w_{q-1}$ be their children in the third row. For $1 \leq l \leq q-1$ let $G_l$ denote the subgraph generated by $w_l$ and all its descendants and let

$$P^l(i) = P^{G_l}(\sigma_{w_l} = i|\sigma_{S(v,2n)} \equiv 1)$$

where $P^{G_l}$ is the model restricted to the subgraph $G_l$ disconnected from the rest of the graph. Note that $P^l(1) = p_{n-1}$. Then

$$p_n = \frac{\sum_{j\in[q+1]^{q-1}} \sum_{k\in[q+1]^{q-1}} \phi(1,j,k) \prod_{l=1}^{q-1} P^l(k_l)}{\sum_{i\in[q+1]} \sum_{j\in[q+1]^{q-1}} \sum_{k\in[q+1]^{q-1}} \phi(i,j,k) \prod_{l=1}^{q-1} P^l(k_l)},$$

where $j = (j_1, \ldots, j_{q-1}) \in [q+1]^{q-1}$ represents the states of $u_1, \ldots, u_{q-1}$, $k = (k_1, \ldots, k_{q-1}) \in [q+1]^{q-1}$ represents the states of $w_1, \ldots, w_{q-1}$ and where $\phi(i,j,k)$ is given by

$$\phi(i,j,k) = \exp\left[\lambda\left(1_{\{i=q+1\}} + \sum_{l=1}^{q-1} 1_{\{j_l=q+1\}}\right)\right.$$
$$\left. - \beta \sum_{l_1<l_2} \eta(j_{l_1}, j_{l_2}) - \beta\left(\sum_{l=1}^{q-1} \eta(j,j_l) + \eta(j_l,k_l) + \eta(j_l,k_{l+1})\right)\right],$$

where $k_q$ is interpreted as $k_1$.

To illustrate $\phi$ let us discuss the case when $\beta = \infty$ where the interactions are hard constraints. Suppose that at most one of $w_1, \ldots, w_{q-1}$ does not have state 1. Then all of the vertices $u_1, \ldots, u_{q-1}$ are adjacent to a state 1 and therefore cannot be in state 1 or in state $q+1$, and since they form a $q-1$-clique they must take every state in $\{2, 3, \ldots, q\}$ exactly once. Since $v$ is connected to $u_1, \ldots, u_{q-1}$ it can not take any of the states $2, 3, \ldots, q$ and since it is adjacent to a 2 it also cannot take the value $q+1$ so it must take state 1.

Applying the same reasoning back in the case of soft constraints with $0 < \beta < \infty$, configurations $(i, j, k)$ where $k$ has at most one state not equal to 1 and where $i \neq 1$ must have at least one $\exp(-\beta)$ in the expansion of $\phi(i, j, k)$ and so $\phi(i, j, k) \leq \exp(q\lambda - \beta)$. The sum of $\prod_{l=1}^{q-1} P^l(k_l)$ over all $k$ with at most one state not equal to 1 is equal to $p_{n-1}^{q-1} + (q-1)p_{n-1}^{q-2}(1 - p_{n-1})$. Finally, for any $(i, j, k)$ we have that $\phi(i, j, k) \leq \exp(q\lambda)$. It follows that for $i \neq 1$

$$\sum_{j \in [q+1]^{q-1}} \sum_{k \in [q+1]^{q-1}} \phi(i, j, k) \prod_{l=1}^{q-1} P^l(k_l)$$
$$\leq e^{q\lambda - \beta}(q+1)^{q-1} + e^{q\lambda}(q+1)^{q-1}(1 - p_{n-1}^{q-1} - (q-1)p_{n-1}^{q-2}(1 - p_{n-1})).$$

On the other hand, if $i = 1$ and $j$ takes every value in $\{2, \ldots, q\}$ exactly once and $k$ is identically 1, then $\phi(i, j, k) = 1$ so

$$\sum_{j \in [q+1]^{q-1}} \sum_{k \in [q+1]^{q-1}} \phi(1, j, k) \prod_{l=1}^{q-1} P^l(k_l) \geq (q-1)! p_{n-1}^{q-1}.$$

Then $p_n \geq f_\beta(p_{n-1})$ where $f_\beta(p)$ is given by

$$(q-1)! p^{q-1}((q-1)! p^{q-1} + e^{q\lambda - \beta} q(q+1)^{q-1}$$
$$+ e^{q\lambda} q(q+1)^{q-1}(1 - p^{q-1} - (q-1)p^{q-2}(1-p)))^{-1}.$$

Now $f_\beta(p)$ converges uniformly on $[0, 1]$ to

$$f(p) = \frac{(q-1)! p^{q-1}}{(q-1)! p^{q-1} + q e^{q\lambda}(q+1)^{q-1}(1 - p^{q-1} - (q-1)p^{q-2}(1-p))}$$

as $\beta$ goes to $\infty$ and $f(1) = 1$ and $f'(1) = 0$. It follows that for some $\varepsilon > 0$, $f(p) > p$ for $p \in (1 - \varepsilon, 1)$. Then for any $p \in (1 - \varepsilon, 1)$ we can find a finite $\beta$ large enough such that $f_\beta(p^*) > p$ for all $p^* \in [p, 1]$ and so

$$\sup_{0 < \beta < \infty} \inf_n P(\sigma_v = 1 | \sigma_{S(v, 2n)} \equiv 1) \geq p$$

which proves (15). But we similarly have

$$\sup_{0 < \beta < \infty} \inf_n P(\sigma_v = 2 | \sigma_{S(v, 2n)} \equiv 2) = 1,$$

which establishes that for large enough $\beta$ there is no unique Gibbs measure.

**7. Remarks.** We will briefly discuss why our example does not immediately provide a counterexample to the conjecture that strong spatial mixing on the $d$-regular tree implies strong spatial mixing on all graphs of maximum degree $d$. Consider the following assignment in the $q$-ary tree with $q \geq 3$. For every vertex $v \in T$ with children $u_1, \ldots, u_q$ set the states of $u_i$ to $i$ for $1 \leq i \leq q-2$. The component of free vertices connected to the root form a binary tree. For large enough $\beta$ because of the conditioning these vertices are likely to take the values $q-1$ or $q$ and are unlikely to take the value $q+1$. Then restricted to the tree of free vertices the distribution is sufficiently close to an antiferromagnetic Ising model which does not have uniqueness for large enough $\beta$. As a result the model does not display strong spatial mixing for sufficiently large $\beta$. We are unable to determine the exact thresholds in $\beta$ for strong spatial mixing on the tree and the graph but do not have any reason to expect that it should produce a counterexample.

**Acknowledgment.** The author would like to thank Elchanan Mossel for suggesting the problem and for his useful comments and advice.

## REFERENCES

[1] BANDYOPADHYAY, A. and GAMARNIK, D. (2006). Counting without sampling: New algorithms for enumeration problems using statistical physics. In *SODA'06: Proceedings of the Seventeenth Annual ACM-SIAM Symposium on Discrete Algorithm* 890–899. ACM, New York.

[2] BRIGHTWELL, G., HÄGGSTRÖM, O. and WINKLER, P. (1999). Nonmonotonic behavior in hard-core and Widom–Rowlinson models. *J. Statist. Phys.* **94** 415–435. MR1675359

[3] DYER, M., SINCLAIR, A., VIGODA, E. and WEITZ, D. (2004). Mixing in time and space for lattice spin systems: A combinatorial view. *Random Structures Algorithms* **24** 461–479. MR2060631

[4] GOLDBERG, L. A., MARTIN, R. and PATERSON, M. (2004). Strong spatial mixing for lattice graphs with fewer colours. In *FOCS'04: Proceedings of the 45th Annual IEEE Symposium on Foundations of Computer Science (FOCS'04)* 562–571.

[5] JUNG, K. and SHAH, D. Inference in binary pair-wise Markov random fields through self-avoiding walks. ArXiv Computer Science e-prints. Available at http://arxiv.org/abs/cs.AI/0610111v2.

[6] MARTINELLI, F. (1999). Lectures on Glauber dynamics for discrete spin models. In *Lectures on Probability Theory and Statistics (Saint-Flour, 1997). Lecture Notes in Mathematics* **1717** 93–191. Springer, Berlin. MR1746301

[7] MOSSEL, E. (2006). Personal communication.

[8] SCOTT, A. D. and SOKAL, A. D. (2005). The repulsive lattice gas, the independent-set polynomial, and the Lovász local lemma. *J. Statist. Phys.* **118** 1151–1261. MR2130890

[9] SOKAL, A. D. (2001). A personal list of unsolved problems concerning lattice gases and antiferromagnetic Potts models. *Markov Process. Related Fields* **7** 21–38. MR1835744

[10] WEITZ, D. (2006). Counting independent sets up to the tree threshold. In *STOC'06: Proceedings of the Thirty-Eighth Annual ACM Symposium on Theory of Computing* 140–149. MR2277139

DEPARTMENT OF STATISTICS
UNIVERSITY OF CALIFORNIA, BERKELEY
BERKELEY, CALIFORNIA 94720
USA
E-MAIL: sly@stat.berkeley.edu